\RequirePackage{fix-cm}

\documentclass{article}

\usepackage{graphicx}
\usepackage[T1]{fontenc}
\usepackage[utf8]{inputenc}
\usepackage[english]{babel}
\usepackage{amsmath}
\usepackage{amsfonts}

\usepackage{amssymb}
\usepackage{latexsym}
\usepackage{bbold}		
\usepackage{url}		
\usepackage{hyperref}		
\usepackage{verbatim}
\usepackage[stable]{footmisc}
\usepackage{framed}
\usepackage{units}
\usepackage{mathtools}		
\usepackage{chngcntr}
\usepackage{rotating}
\usepackage[table]{xcolor}
\usepackage{comment}		
\usepackage{algorithm,algpseudocode}
\definecolor{grey}{gray}{0.9}

\DeclareMathOperator{\R}{\mathbb{R}}					
\DeclareMathOperator{\1}{\mathbb{1}}					
\DeclareMathOperator{\argmax}{argmax} 					
\DeclareMathOperator{\argmin}{argmin} 					
\DeclareMathOperator{\vect}{vec} 					
\DeclareMathOperator{\orth}{orth} 					

\renewcommand\Re{\operatorname{Re}}
\usepackage{draftwatermark}
\SetWatermarkText{Preprint}
\SetWatermarkScale{1} 

\begin{document}

\title{Data-Driven Combined State and Parameter Reduction for Extreme-Scale Inverse Problems}


\author{Christian Himpe\thanks{Contact: \href{mailto:christian.himpe@uni-muenster.de.de}{\nolinkurl{christian.himpe@uni-muenster.de}}, \href{mailto:mario.ohlberger@uni-muenster.de}{\nolinkurl{mario.ohlberger@uni-muenster.de}}, Institute for Computational and Applied Mathematics at the University of M\"unster, Einsteinstrasse 62, D-48149 M\"unster, Germany} \and Mario Ohlberger\footnotemark[1]{}}


\date{}

\maketitle

\begin{abstract}
In this contribution we present an accelerated optimization-based approach for combined state and parameter reduction of a parametrized linear control system which is then used as a surrogate model in a Bayesian inverse setting. 
Following the basic ideas presented in [Lieberman, Willcox, Ghattas. Parameter and state model reduction for large-scale statistical inverse settings, SIAM J. Sci. Comput., 32(5):2523-2542, 2010], our approach is based on a generalized data-driven optimization functional in the construction process of the surrogate model and the usage of a trust-region-type solution strategy that results in an additional speed-up of the overall method.
In principal, the model reduction procedure is based on the offline construction of appropriate low-dimensional state and parameter spaces and an online inversion step based on the resulting surrogate model that is obtained through projection of the underlying control system onto the reduced spaces. 
The generalization and enhancements presented in this work are shown to decrease overall computational time and increase accuracy of the reduced order model and thus allow an application to extreme-scale problems. 
Numerical experiments for a generic model and a fMRI connectivity model are presented in order to compare the computational efficiency of our improved method with the original approach.

\end{abstract}

\section{Introduction}
Many physical, chemical, technical, environmental, or bio-medical applications require the solution of inverse problems for parameter estimation and identification. 
This is in particular the case for complex dynamical systems where only experimental functional data is accessible via measurements. 
In neurosciences, a particular application is e.g. the extraction of effective connectivity in neural networks from measured data, such as data from electroencephalography (EEG) or functional magnetic resonance imaging (fMRI).

If a network with many states is considered, the corresponding large-scale inverse problem is often only accessible in reasonable computational time if model reduction is applied to the underlying control system. 
Moreover, the measured data is subject to statistic errors such that it is reasonable to apply a Bayesian inference approach which tries to identify a distribution on the underlying parameters, rather than computing deterministic parameter values. 
In the context of connectivity analysis in neurosciences such an inversion approach has been established by Friston and his collaborators in recent years under the synonym Dynamic-Causal-Modeling (DCM) \cite{friston03}. 
We in particular refer to \cite{stephan07} and the references therein for further reading. 

In this contribution we will consider model reduction for Bayesian inversion of general linear control systems which in particular includes the above mentioned application scenario in neurosciences, but also other input-output systems that may be obtained, for example, in a partial differential equations setting after discretization. 

Control systems usually comprise some internal state ($x$), allow external input ($u$), transform the input and state to observed output ($y$) and provide different configurations through parameters ($\theta$).
Generally, a control system therefore consists of a dynamical system with vector field $f$ and an output functional $g$.
In the finite dimensional case (e.g. after discretization) we may consider $x(t) \in R^N$ for $t \in \tau := [0,T]$, and a control system given by:
\begin{align*}
 \dot{x}(t) &= f(x(t),u(t),\theta), \qquad  \forall t \in \tau,\\
       y(t) &= g(x(t),u(t),\theta),\qquad  \forall t \in \tau,
\end{align*}
with external input $u(t) \in R^J$, output $y(t) \in R^O$ and parameters $\theta \in R^P$. 
Naturally, the system is equipped with an initial condition for the state, i.e. $x(0) = x_0 \in \R^N$.

In this contribution, as a general underlying model, linear control systems of the following form are considered:
\begin{align}\label{eq:lti}
 \dot{x} &= A(\theta) x + B u,\\
       y &= C x,\notag
\end{align}
where now, we have neglected the dependence on time for the ease of exposition. 
Here, $A(\theta) \in \R^{N \times N}$ denotes the parametrized system matrix, $B \in \R^{N \times J}$ an input matrix , and $C \in \R^{O \times N}$ an output matrix.
The matrix $A(\theta)$ is fully parametrized, meaning each component is an individual parameter, thus $\theta \in \R^{N^2}$ is mapped to the system matrix by:
\begin{align*}
 {\vect^{-1}}:\R^{N^2} &\to \R^{N \times N},\\
 \theta_i &\mapsto A_{(i/N,i\%N)}.
\end{align*}
These $N^2$ parameters are unknown up to their prior distribution at the time of reduction, and for models with $N >\!\!> 1$ the number of parameters might make an estimation prohibitively computationally expensive.

For bounded input, the systems stability is characterized by the eigenvalues of $A(\theta)$.
Hence, the parameters composing $A$ need to be estimated in a manner that the system remains stable, which can be considered an additional prior information that needs to be taken into account during the (Bayesian) inversion.

Reducing the parameter space, e.g. the number of independent connections in a brain connectivity model, before the parameter estimation, can lower the complexity of optimization significantly.
Alternatively, reducing the state space, in our example the number of brain regions or nodes, cuts the computational cost of the necessary integrations during the optimization.
The resulting reduced order model approximates the original model up to an error introduced by the reduction procedure.
In this contribution, a {combined reduction} of state and parameter space is presented, allowing the swift inversion of large-scale and even extreme-scale parametrized models.

The inversion procedure incorporating model reduction is usually arranged in two steps.
In a first step, called the {offline} phase, the underlying parametrized model is reduced; here in states and parameters.
Second, in an {online} phase, the reduced models parameters are estimated to fit the observed experimental data.

In Bayesian inverse problems \cite{stuart10,biegler11} we aim at estimating a probability distribution of the unknown parameters instead of the parameter values directly.
If we denote the given data by $y_d$, the inverse approach computes the so called posterior distribution, $P(\theta | y_d)$, i.e. the probability distribution 
of the parameters under the given data. The basic underlying mathematics to achieve this goal is Bayes' rule:
\begin{align*}
 P(\theta | y_d) = \frac{P(y_d | \theta)P(\theta)}{P(y_d)};
\end{align*}
with the prior distribution $P(\theta)$, reflecting the preliminary knowledge about the parameter distribution.
$P(y_d|\theta)$ is the likelihood that we observe $y_d$, given a parameter vector $\theta$. 
The likelihood is estimated in the online phase and requires forward integrations of the underlying dynamical system. 
$P(y_d)$ is the model evidence - the probability of the data.

Up to the normalizing factor $P(y_d)$, the posterior distribution is proportional to $P(y_d | \theta)P(\theta)$, i.e.
\begin{align*}
 P(\theta | y_d) \propto P(y_d | \theta)P(\theta).
\end{align*}

In our problem setting, the observed data $y_d$ is presumed to contain some additive noise $\epsilon$, i.e.
\begin{align*}
 y(\theta) + \epsilon \approx y_d.
\end{align*}
In the scope of this work Gaussian noise $\epsilon = N(0,v)$ is assumed. 
Thus, in this fully Gaussian setting, all probability distributions can be specified in terms of mean and covariance.
Therefore, given a prior Gaussian distribution $P(\theta)~=~N(\kappa,K)$ and some experimental data $y_d$, a distribution proportional to the posterior distribution can be computed by the Gaussian likelihood $P(y_d|\theta)~=~N(\lambda,\Lambda)$ and prior:
\begin{align*}
 P(\theta | y_d) \propto \exp(-\frac{1}{2}\|y(\theta) - y_d\|^2_{\Lambda^{-1}} - \frac{1}{2}\|\theta - \kappa\|^2_{K^{-1}}).
\end{align*}
Then, the Maximum-A-Posteriori (MAP) estimator is given by:
\begin{align*}
 \theta_{MAP} &= \argmax \exp(-\frac{1}{2}\|f(\theta) - y_d\|^2_{\Lambda^{-1}} - \frac{1}{2}\|\theta - \kappa\|^2_{K^{-1}}) \\ 
              &= \argmin \left( \frac{1}{2}\|f(\theta) - y_d\|^2_{\Lambda^{-1}} + \frac{1}{2}\|\theta + \kappa\|^2_{K^{-1}}\right)
\end{align*}
which is computed in the online-phase after we have constructed a suitable reduced order model.

Model reduction of parametrized control systems has been investigated in recent years through various approaches. 
Early interpolation ideas were presented in \cite{WMGGD99} as well as moment matching techniques e.g. in \cite{FB07}.
Some recent approaches comprise interpolatory schemes with sparse grids \cite{BB08}, superposition of locally reduced models \cite{LE09} or matrix interpolation \cite{eid11,amsallem11}. 
Another approach comes from the reduced basis techniques. 
Here, in particular the POD-Greedy algorithm has been established \cite{haasdonk08,haasdonk13} in the context of parametrized partial differential equations and transfered to dynamical systems in \cite{haasdonk11}. 
We also refer to \cite{nguyen11} for an application of reduced basis model reduction in an Bayesian inverse setting.

While all the above mentioned approaches concentrate on state space reduction only, there are also very recent approaches towards simultaneous reduction of state and parameter spaces in the context of large-scale inverse problems. 
Concerning gramian-based combined parameter and state reduction we in particular refer to \cite{himpe13} and the references therein. 
In this contribution, however, we are concerned with optimization-based combined state and parameter reduction.
Our approach is mainly based on \cite{lieberman10} where a procedure to concurrently reduce state and parameter spaces of linear control systems in a Bayesian inversion setting has been introduced.
The iterative improvement of a projection as used in \cite{lieberman10} is also used in \cite{antoulas09} from a gramian-based perspective.
This method is related to the state reduction from \cite{buithanh07}, the parameter reduction from \cite{buithanh08} and is generally related to the Hessian-based model reduction ansatz as described in \cite{bashir08}. 
We also refer to \cite{lieberman12} for a corresponding goal-oriented optimization based approach.

Starting from the ideas in \cite{lieberman10,lieberman12} we base our approach in this contribution on a generalized data-driven optimization functional and present enhancements of the reduction procedure using ideas from \cite{haasdonk08,arian00}.

The article is organized as follows.
In section~\ref{combined}, a short review of the model reduction procedure from \cite{lieberman10} is given.
In section~\ref{enhan}, we first take into account a data-driven optimization functional, and then present enhancements to the existing approach that result in a reduction of offline computational time.
Section~\ref{impl} summarizes the implementation of the presented algorithm and its extensions.
Finally, we evaluate the resulting model reduction approach in numerical experiments.
The methods are tested and compared on a generic model and an fMRI connectivity model for synthetic neuronal activity in section~\ref{numer}.

\section{Combined State and Parameter Reduction}\label{combined}
In this section, the method from \cite{lieberman10} is briefly reviewed and annotated.
For large-scale inverse problems, model reduction becomes relevant to ensure reasonable optimization durations.
Commonly, Galerkin or Petrov-Galerkin projections are employed for the low-rank projections of state \cite{bashir08} or parameter \cite{buithanh08} spaces.
The simultaneous reduction of state and parameter space is based on Galerkin projections $V$, with:
\begin{align*}
 V^T V = \1.
\end{align*}

The reduced model is of lower order $M \ll N$ than the original full-order model.
For the reduced states $x_r \in \R^M$ the associated low-rank control system is derived from the original models components:
\begin{align*}
 \dot{x}_r &= A_r(\theta_r) x_r + B_r u, \\
 y_r &= C_r x_r,
\end{align*}
with a reduced initial condition $x_{r,0} = V^T x_0$ and the reduced components:
\begin{align*}
 \theta_r &= P\theta \in \R^{M^2}, \\
 A_r(\theta_r) &= V^TA(P^T\theta_r)V \in \R^{M \times M}, \\
 B_r &= V^TB  \in \R^{M \times J}, \\
 C_r &= CV  \in \R^{O \times M}.
\end{align*}
The required state projection $V \in \R^{M \times N}$ and parameter projection $P \in \R^{M^2 \times N^2}$ are determined iteratively.
This iterative assembly of the projection matrices is based on a greedy algorithm \cite{lieberman07,buithanh08}, that optimizes the error between the high-fidelity original and the low-dimensional reduced model.
In each iteration a set of reduced parameters is determined by maximizing the error (see \cite{bashir08}) between the original and the reduced models output using the following objective function
\begin{align*}\label{opt_simple}
 J(\theta) &= \alpha \Delta(y_r(\theta_r)) - \beta \| \theta_r \|_{K^{-1}}^2
\end{align*}
with suitable weights $\alpha, \beta \in [0,1]$.
Here $\Delta(y_r(\theta_r))$ denotes a measure for the error in the output between the reduced model evaluated at the reduced parameter and the full underlying model with the high dimensional parameter. 
In the original approach \cite{lieberman10} the output error measure is chosen as $\Delta(y_r(\theta_r)) = \| y_r(\theta_r) - y(\theta) \|_2^2$. 
The regularization in the second term utilizes the prior covariance matrix $K \in \R^{N^2 \times N^2}$.
In case the covariance matrix $K$ is diagonal, the inverse of the covariance matrix, the precision matrix, $K^{-1}$ can be computed by inverting each diagonal component of $K$. 
Thus, the second summand simply regularizes this maximizer in terms of the provided prior distribution:
\begin{align*}
 \|\theta_r\|_{K^{-1}}^2 = (P \theta_r)^T K^{-1} (P \theta_r),
\end{align*}
which is a weighted $2$-norm.
This type of regularization makes use of the prior covariance and thus penalizes parameters of low probability with respect to the prior information.

The presented model order reduction method relies on optimization to fit the reduced parameters optimally.
As in \cite{lieberman10} and \cite{bashir08} a greedy method can be employed.
Yet, one is not restricted to the output error in the $2$-norm. 
Alternatively, the least-absolute-deviation method ($1$-norm minimization), or the least-maximum-deviation ($\infty$-norm minimization) can be used, i.e. 
\begin{align*}
 \Delta(y_r(\theta_r)) = \| y_r(\theta_r) - y(\theta) \|_1 \quad {\rm or} \quad 
 \Delta(y_r(\theta_r)) = \| y_r(\theta_r) - y(\theta) \|_\infty.
\end{align*}
More generally, each $p$-norm can be used inside the objective function, which results in a generalized objective functional with 
$
 \Delta(y_r(\theta_r)) = \| y_r(\theta_r) - y(\theta) \|_p^p
$
for some $p\in [1, \infty]$. 
The regularization term remains unchanged, since it is solely based on the prior distribution.
However, since the evaluation of the output error in the $p$-norm requires high dimensional solves in each step of the optimization loop, it is in general advisable 
to replace the true output error by some a posteriori error estimator, as e.g. derived and suggested in \cite{haasdonk08,grepl11}.
Such an a-posteriori error estimator can be computed more efficiently since it does not require full-order time integrations.

Each iteration requires computing the reduced model, based on the last iterations' projection matrices $\{V, P\}$, as well as the greedy optimization of $J$ based on the integration of the full and the reduced model:
\begin{align*}
 \theta_I &= \argmax J(\theta) \\ 
          &= \argmax \;\;\: \Delta(y_r(\theta_r)) - \beta \| \theta_r \|_{K^{-1}}^2 \\
          &= \argmin \; -\Delta(y_r(\theta_r)) + \beta \| \theta_r \|_{K^{-1}}^2 .
\end{align*}
The resulting reduced parameters constitute the next basis vector being orthogonalized\footnote{The orthogonalization  of state and parameter projection can be accomplished by various algorithms for example Gram-Schmidt, Householder-Reflections, Givens-Rotation or Singular Value Decomposition} into the parameter projection $P$:
\begin{align*}
 \begin{pmatrix} P_I, & \theta_I \end{pmatrix} &= QR \\ \Rightarrow P_{I+1} &= Q.
\end{align*}

As the dynamic system is linear and time-invariant, a simulation of such system is equivalent to solving a system of linear equations \cite{bashir08}:
\begin{align}
 \begin{pmatrix} \1 & & & \\ A(\theta) & \1 & & \\ & \ddots & \ddots & \\ & & A(\theta) & \1 \end{pmatrix}  
 \begin{pmatrix} x_0 \\ x_1 \\ \vdots \\ x_T \end{pmatrix} = 
 \begin{pmatrix} B u(0) \\ B u(1) \\ \vdots \\ B u(T) \end{pmatrix}.
\end{align}
Some selection $\overline{x(\theta_I)}$ of the solution time series $x(\theta_I)$ is then incorporated into the state projection by orthogonalization:
\begin{align*}
 \begin{pmatrix} V_I, & \overline{x(\theta_I)} \end{pmatrix} &= QR \\ \Rightarrow V_{I+1} &= Q.
\end{align*}
In \cite{bashir08} POD-modes are selected as $\overline{x(\theta_I)}$ to be included into the projection.
A more simple but numerically very efficient approach would be the mean of the time series, as used in \cite{lall99} and \cite{lall02}:
\begin{align*}
 \overline{x} = \frac{1}{T} \int_0^T x(t) dt,
\end{align*}
and in the discrete case:
\begin{align*}
 \overline{x} = \frac{\Delta t}{T} \sum_{i=1}^T x_i.
\end{align*}
In this contribution, we suggest to select $\overline{x(\theta_I)}$ from a truncated POD of the orthogonal projection error 
in the time series $x(\theta_I)$ with respect to the reduced state space of the preceding iteration. This approach follows the idea of the 
POD-Greedy procedure proposed in \cite{haasdonk08}.

Using the projections $P$ and $V$, the next iteration is performed.

The parameter projection $P$ can be initialized with a constant vector as described in \cite{lieberman10}, yet a more natural choice is the prior mean, assuming it is not identical to zero, $\theta_{\text{prior}} \not\equiv 0$.
A prior information that is usually implicitly assumed is the underlying systems stability.
Hence, without any other prior information one could at least choose $\theta_{\text{prior}}=N(-\mathbb{1}_N,\1_{N \times N})$ as uninformative priors, suggesting stability:
\begin{align*}
 P_0 = \theta_0 = \vect(-\1_{N \times N}).
\end{align*}
From this initial choice for the parameters the full order system is sampled and the state projection $V$ is initialized, for example, by the mean over time of the states:
\begin{align*}
 V_0 = \overline{x(\theta_{\text{prior}})}.
\end{align*}

In summary, the complete reduction algorithm is given by the pseudo-code listing of algorithm~\ref{alg:a1}.
\begin{algorithm}
 \caption{Original Combined State and Parameter Reduction}
 \label{alg:a1}
 \begin{algorithmic}
   \State $\theta_0 \leftarrow \theta_{\text{prior}}$
   \State $P_0 \leftarrow \theta_0$
   \State $V_0 \leftarrow \overline{x(\theta_0)}$
   \For{I=1:R}
    \State $A_r(\theta_r) \leftarrow V_I A(P_I \theta_r) V_I$
    \State $B_r \leftarrow V_I B$
    \State $C_r \leftarrow CV_I$
    \State $\theta_I \leftarrow \argmax J(\theta_{I-1})$
    \State $P_{I+1} \leftarrow \orth(P_I,\theta_I)$
    \State $V_{I+1} \leftarrow \orth(V_I,\overline{x(\theta_I)})$
   \EndFor
 \end{algorithmic}
\end{algorithm}

In algorithm \ref{alg:a1}, $x(\theta_I)$ describes a snapshot of the states for parameters $\theta_I$.
$\overline{x(\theta_I)}$ corresponds to the selection from the states over time and the \texttt{orth} method orthogonalizes a given matrix.
The \texttt{argmax} method represents an optimization procedure, which in the scope of this work is given by, but is not restricted to, an unconstrained optimization.

A reconstruction, after the inference, of the parameters, using the above described reduction procedure of states and parameters, is also accomplished by the computed projections $V$ and $P$.
Due to the use of Galerkin projections, the inverse projection of $P$ and $V$ is given by their transpose $P^{-1} = P^T$, $V^{-1} = V^T$, thus:
\begin{align*}
 \theta_{r,\text{prior}} &= P \theta_{\text{prior}} \approx P^T \theta_{\text{prior}}, \\
 S_{r,\text{prior}} &= V^T S_{\text{prior}} V \approx V S_{\text{prior}} V^T.
\end{align*}

With this \textbf{combined reduction} of parameter and state space, using the parameter projection $P$ and the state projection $V$, improves the inversion procedure not only by shortening the integration durations due to the state reduction, but also by decreasing the number of optimizable parameters.

\section{Incorporating Data-Misfit and Trust-Region Methods}\label{enhan} 
In this section we first include experimental data into the model reduction procedure.
Second, a trust-region like approach is presented to shorten the duration of the reduced order model construction.
These are aimed to accelerate the assembly of the model reducing projections $P$, $V$ and thus shorten the overall offline phase.

\subsection{Incorporating Data-Misfit}\label{data}
Inverse problems require some experimental data $y_d$ to which a models parameters are fit during the inversion in the online phase.
We thus incorporate the data-misfit $\|y_r - y_d\|$ into the objective functional in the combined reduction procedure during the offline phase. 
This is motivated by two arguments. 
While in the original approach, the surrogate model is constructed to be accurate within the whole parameter space, it can now be tailored towards the parameter ranges that are related to the solution of the inverse problem, which might result in lower dimensional surrogate models that lead to more accurate results in the inversion.
If only data-misfit is considered, the objective functional would read (cf. \cite{lieberman12})

\begin{align*}
 J_d(\theta) = -\| y_r(\theta_r)-y_d \|,
\end{align*}
where $\| \cdot \|$ denotes a suitably chosen norm. 
Note, that the usage of the data-misfit will save computational time during the offline phase, as no full order integrations of the model is needed to evaluate the objective functional.

The resulting optimization problem will most likely be ill-posed. Therefore, we will keep the regularization from eq.~\ref{opt_simple}. 
Furthermore, for the solution of the inverse problem it is advisable to prepare the reduced model in the offline phase in such a way that also parameter ranges in the neighborhood of the optimal parameters associated with the measured data are taken into account. 
Our suggested generalized objective functional $\hat{J}_p$ is a combination of $J_p$ with the data-misfit $J_d$. 
Using weighting parameters $\alpha, \beta, \gamma$; we thus obtain
\begin{align*}
 \hat{J}_p(\theta) = \alpha \Delta_p(y_r(\theta_r)) - \beta \|\theta_r\|_{K^{-1}}^2 - \gamma \| y_r(\theta_r)-y_d \|_p^p.
\end{align*}
A possible choice that equilibrates the influence of all three terms would be $\alpha=\beta=\gamma=1/3$.
As a result, the reduced models parameters are determined to encourage matching the reduced model to the provided experimental data instead to a general reduced model with the given priors.
To further enforce the fitting of the reduced model to the measured output, we could choose $\beta = \gamma = 1 - \alpha$ with $\alpha$ very small or even $\alpha=0$ which results in an objective functional without the main maximization term of eq.~\ref{opt_simple}, i.e.
\begin{align*}
 \tilde{J}_p(\theta) = - \beta \|\theta_r\|_{K^{-1}}^2 - \gamma \| y_r(\theta_r)-y_d \|_p^p.
\end{align*}
Apart from an expected higher accuracy in matching the experimental data, which is the ultimate goal in a later online phase of the inversion procedure, this massively lowers the computational load, since no sampling of the full order model is required.
Yet, due to the usage of specific experimental data in the reduction procedure, the resulting reduced model is only valid for fitting this particular data.

\subsection{Trust-Region Strategy}\label{trust}
The original model reduction algorithm keeps the dimension of the reduced parameter vector fixed and iterates until convergence or predetermined reduced order.
Yet the dimension of the to-be-estimated parameter vector largely determines the offline duration, due to the required integrations during optimization in each iteration.
This can be counteracted by a trust-region like strategy, which is loosely related to \cite{arian00} and \cite{kelley99}, where trust-region methods are applied to POD based model reduction.
In this contribution the basic trust-region algorithm (\cite[Ch. 6]{conn00}) is simplified to allow a swift computation by removing the acceptance step.
Due to the optimization during each iteration of the reduction process, which itself iterates until some acceptable bound is reached, an extra acceptance step is not required. 
Since the dimension of the parameter vector $\tilde{\theta}$ varies over iterations, so for the $I$th iteration:
\begin{align*}
 \dim(\tilde{\theta}_I) = I,
\end{align*}
and an additional mapping from the trust-region parameter vector $\tilde{\theta}_I$ to the full-order parameter vector $\theta_I$ is required.
This mapping enables the orthogonalization of the current iterations parameter vector into the parameter projection.
A simple mapping is given by:
\begin{align*}
 \varphi:\R^I &\to \R^{N^2}, \\
 \tilde{\theta}_I &\mapsto \begin{pmatrix} \tilde{\theta}_I \\ 0 \\ \vdots \\ 0 \end{pmatrix} = \theta_I.
\end{align*}
For the next iteration $I+1$, the parameter vector is extended by:
\begin{align*}
 \tilde{\theta}_{I+1} = \begin{pmatrix} \tilde{\theta}_I \\ 0 \end{pmatrix} \in \R^{\dim(\tilde{\theta}_I)+1}.
\end{align*}
The trust-region radius is initially set to dimension one, thus instead of initializing the parameter vector with a constant vector or the prior means, it is set to (scalar) one:
\begin{align*}
 \tilde{\theta}_0 &= 1 \ \in \R.
\end{align*}
Yet, the first column of the parameter projection is still initially set to the prior means:
\begin{align*}
 P_0 &= (\overline{x(\theta_{\text{prior}})})^T \in \R^{N^2}.
\end{align*}
Then, in the first iteration of the reduction the (scalar) parameter $\theta_0$ is computed which approximates the full system $A(\theta)$ best,
\begin{align*}
 \tilde{\theta}_1 &= \argmax(J(\theta)) \in \R  \\ 
 \Rightarrow P_1 &= \orth \begin{pmatrix}P_0, & \varphi(\tilde{\theta}_1)\end{pmatrix}.
\end{align*}
In each subsequent iteration an additional dimension of the parameter space is added to the trust-region radius.
For example, the second iteration optimizes two parameters, the third optimizes three, etc until the given reduced parameter space dimension is reached.

At each instance the full parameter vector is required, it is projected by the inverse (which equates to the transposed, due to the orthogonalization) of the parameter projection $P$:
\begin{align*}
 \theta \approx P^{-1} \varphi(\tilde{\theta_r}) = P^T \varphi(\tilde{\theta_r}).
\end{align*}
Hence, starting with a single scalar parameter, that is initialized with $1$, the reduced parameter vector is assembled by iteratively incrementing the dimension.

With this enhancement, in each iteration an optimization problem of lesser or (once only) equal dimension than the original algorithm.
Due to the smaller size of the optimization problem the offline time is massively lowered.

\section{Implementation}\label{impl}
The trust-region enhancement, together with the data-misfit enhancement from section~\ref{data} are modularly included into the algorithm; listing \ref{alg:a2} showcases the new algorithm.

\begin{algorithm}
 \caption{Enhanced Combined State and Parameter Reduction}
 \label{alg:a2}
 \begin{algorithmic}
   \If{trust-region}
   \State $\theta_0 \leftarrow 1$
   \Else
   \State $\theta_0 \leftarrow \theta_{\text{prior}}$
   \EndIf
   \State $\overline{x(\theta_0)} \leftarrow \begin{cases} \text{pod}(x(\theta_0)) \\ \text{mean}(x(\theta_0)) \\ \text{pod-greedy}(x(\theta_0)) \end{cases}$ 

   \State $P_0 \leftarrow \theta_{\text{prior}}$
   \State $V_0 \leftarrow \overline{x(\theta_0)}$
   \For{I=1:R}
    \State $A_r(\theta_r) \leftarrow V_I A(P_I \theta_r) V_I$
    \State $B_r \leftarrow V_I B$
    \State $C_r \leftarrow CV_I$
    \State $J \leftarrow \begin{cases} J_{\text{original}}(\alpha,\beta;\theta) \\ J_{\text{data-driven}}(\alpha,\beta,\gamma;\theta) \\ J_{\text{data-only}}(\beta,\gamma;\theta) \end{cases}$
    \State $\theta_I \leftarrow \argmax J(\theta_{I-1})$
    \If{trust-region}
    \State $\theta_I \leftarrow \varphi(\theta_I)$
    \EndIf
    \State $\overline{x(\theta_I)} \leftarrow \begin{cases} \text{pod}(x(\theta_I)) \\ \text{mean}(x(\theta_I)) \\ \text{pod-greedy}(x(\theta_I)) \end{cases}$
    \State $P_{I+1} \leftarrow \orth(P_I,\theta_I)$
    \State $V_{I+1} \leftarrow \orth(V_I,\overline{x(\theta_I)})$
   \EndFor
 \end{algorithmic}
\end{algorithm}

The algorithm from code listing \ref{alg:a2} is implemented under the name \texttt{optmor} - \textbf{opt}imization-based \textbf{m}odel \textbf{o}rder \textbf{r}eduction.
The source code is available from: \url{http://j.mp/optmor} under an open-source license and is compatible with OCTAVE and MATLAB.
For compatibility reasons the estimation algorithm employed during the reduction is an unconstrained optimization, but can easily be replaced with a constrained optimization function (\cite{buithanh08}, \cite{lieberman10}).
To remain configurable the here described enhancements can each be used optionally; either individually or in combination.
Additionally, the usage of a source term $F$ and a feed-forward matrix $D$ is implemented as well, allowing models of the form:
\begin{align*}
 \dot{x} &= A(\theta) x + B u + F, \\
 y &= C x + D u,
\end{align*}
to be reduced.

The interface of the \texttt{optmor} program is given by:
\begin{align*}
 \text{\texttt{\{P,V\} = optmor(p,A,B,C,D,F,T,R,X,U,S,q,y)}};
\end{align*}

with $p$ being a vector containing the parameter prior mean.
The argument $A$ is a function handle to a mapping from a given parameter vector to a system matrix using the signature: \texttt{A = @(p)}.
If $A = 1$ the inverse mathematical vectorization map, $\vect^{-1}$, will be assumed as parameter mapping.
$B$ is the input matrix, $C$ is the output matrix, $D$ is the feed-forward matrix and $F$ is the source term.
The vector $T$ is a three component vector holding the start time, time step and end time, while the scalar $R$ holds the targeted reduced dimension.
Furthermore, the vector $X$ represents the initial value and the matrix $U$ provides the input or control, for each time step.
The argument $S$ holds the associated prior covariance matrix; for $S = 1$ a unit covariance matrix is assumed.
$q$ is a six component vector holding the configurable options.
Optionally, $y$ may hold experimental output time series, required for the data-misfit enhancements.
The algorithm returns two projection matrices, $P$ and $V$, for parameter and state projections respectively.

For selecting $\overline{x(\theta)}$ from a snapshot also the POD-Greedy method \cite{haasdonk08} using the error system is implemented.
As described above, in \cite{bashir08} a system of linear equations is solved to simulate the forward model.
Alternatively, a single-step solver, like a Runge-Kutta, or multi-step solver, like Adams-Bashforth, can be utilized to solve the system.
An advantage of using such solvers is the lesser memory requirements opposed to solving a linear system with dimension being the product of states and time-steps.

\section{Numerical Results}\label{numer}
To demonstrate the capabilities of this approach in combined state and parameter reduction, two types of models are tested.
First, a generic control system as described in eq.~\ref{eq:lti} is tested.
Second, the combined reduction is applied to a linearized system for the inversion of fMRI data to deduce connectivity between brain regions \cite{kamrani12,stephan07}.
Lastly, an extreme-scale problem is tested as well as an evaluation of the effectivity of the reduction method for different configurations.

\subsection{Online Phase}
In the online phase, the estimation of the parameter (distribution) is accomplished by a least-squares minimization of the residual between reduced order model output and experimental data.
The objective function employed in the optimization of the full-order model is given by:
\begin{align*}
 J_{\text{\small online}}(\theta) &= \| y(\theta)-y_d \|_2^2 + \|\theta\|_{K^{-1}}^2,
\end{align*}
whereas for the reduced models an adapted objective function of the following form is utilized:
\begin{align*}
 \tilde{J}_{\text{\small online}}(\theta_r) &= \| y_r(\theta_r)-y_d \|_2^2 + \|P^T\theta_r\|_{K^{-1}}^2.
\end{align*}
Here, the parameter estimation is performed with an unconstrained (least-squares) optimization with regularization for the full-order and reduced order models.

\subsection{Generic Linear Control System}
The \texttt{optmor} implementation is tested with a generic linear control system.
As mentioned above, we will assume $A(\theta)$ is fully parametrized, hence $\theta \in \R^{N^2}$.
The number of states is varied with $N=\{ 9,16,25,36\}$, and thus $P = N^2$; while the number of inputs and outputs is fixed to $I = O = \sqrt{N}$.
Systems with these dimensions are generated randomly, but with ensured stability\footnote{For the real part of the eigenvalues of the system matrix holds: $\Re(\lambda(A(\theta))) < 0$} for each set of experiments.
Input $u(t)$ will be given by a simple delta impulse, while the source term $F$ and the feed-forward matrix $D$ are set to zero.
The prior mean of the parameters is set to $-1$ on the diagonal, $0$ off the diagonal of $A(\theta) \Rightarrow \theta_{\text{prior}} = \vect(-\1)$, thus ensuring initial stability of the system; while the prior covariance is set to the unit matrix $S = \1$.
For all the reductions the number of iterations is fixed to the number of outputs, which allows comparison of the different combined reduction variants of the same reduced order.
For the online phase and the data-driven extensions, simulations using the original parameters with added Gaussian noise are used.
The source code, used to conduct the experiments below, can be found at \url{http://j.mp/optmor}.

As a baseline, the full-order model is estimated without employing any reduction.
Since in this case the full-order models (high-dimensional) parameters are approximated, there is only an online phase.

\begin{figure}[h!]
 \includegraphics[scale=0.38]{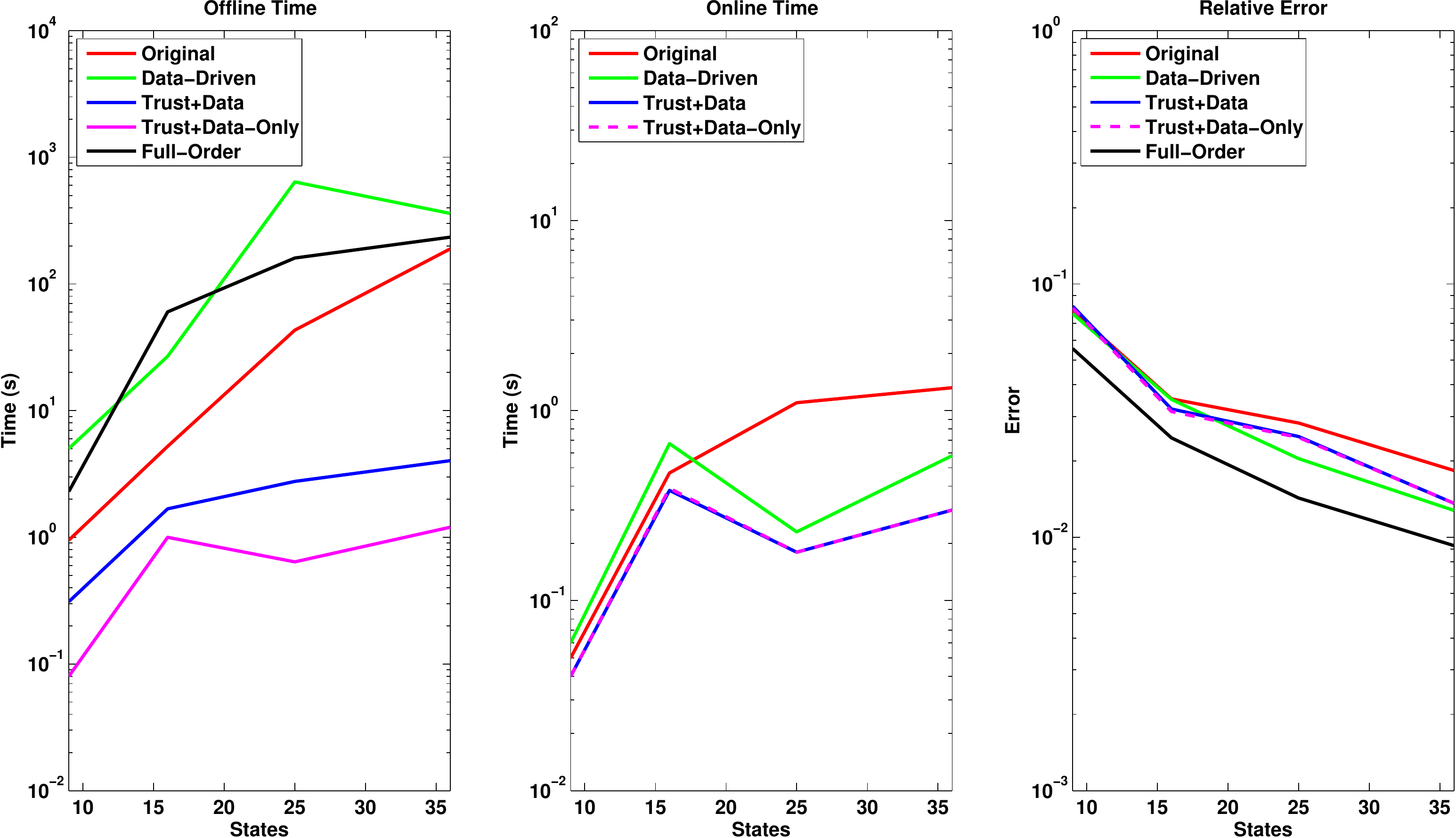}
 \caption{Offline time, online time and relative error for the original, data-driven, trust-region and combined data-driven/trust-region reductions.}
 \label{fig:trust}
\end{figure}

Next, the presented data-driven (section~\ref{data}) and trust-region (section~\ref{trust}) extensions are tested and compared individually and in combination for this parametrized linear control system.
In figure~\ref{fig:trust} the offline and online durations as well as the relative error in outputs is shown for the enhancements in comparison to the full-order optimization and the original reduction method.

The additional data-misfit term, of the data-driven enhancement, increases the offline phase duration, but reduces the relative error and online time.
The newly introduced trust-region enhancement greatly reduces offline phase duration as predicted and also to a lesser degree the online time compared to the original and data-misfit method.
Compared to the original methods' \cite{lieberman10} relative error, the relative error behaves slightly better.
Combining the data-driven and trust-region approach significantly shortens the offline phase, yielding a slightly higher relative error compared to the data-driven method but below the original algorithms error.
Lastly, excluding the minimization term from the objective function ($\alpha=0$), again reduces the offline time while not affecting online time and relative error.
Thus, especially the combination of the data-driven and trust-region enhancements massively accelerates the reduction process.

Assessing the effectivity of the reduction, using the combined data-driven and trust-region approach, in terms of total (offline and online) time compared to the full-order solution and the original algorithm results in speed-ups of up to two orders of magnitude as listed in table~\ref{eff1}.

\begin{table}\centering
\begin{tabular}{c|c|c}
 State Dimension & Speed Up (vs. Full-Order) & Speed Up (vs. Original) \\ \hline
  9 &  $19.71$ &   $8.33$ \\
 16 &  $43.31$ &   $4.09$ \\
 25 & $195.17$ &  $54.04$ \\
 36 & $156.22$ & $127.00$ 
\end{tabular}
\caption{Effectivity of the exclusive combined data-driven and trust-region method for the generic model.}
\label{eff1}
\end{table}

\subsection{fMRI Connectivity Model}
A method to infer connectivity for different regions of the brain based on experimental data recorded by fMRI or fNIRS is known as Effective Connectivity \cite{stephan07}.
There are two sub-models composing the underlying model of Effective Connectivity which is a concept that is closely related to Dynamic-Causal-Modeling \cite{friston03}.
The dynamic sub-model represents the network of the observed brain regions by a controlled linear system:
\begin{align*}
 \dot{x} = A(\theta_{\text{dyn}})_{\text{dyn}} x + B_{\text{dyn}} u
\end{align*}

The forward sub-model converts each state $x_i$ of the dynamic sub-model to the observed measurements.
In the case of fMRI observation the forward-submodel is given by the nonlinear hemodynamic model \cite{friston03}:
\begin{align*}
 \dot{s_i} &= x_i - \frac{1}{\tau_s}s - \frac{1}{\tau_f}(f_i-1), \\
 \dot{f_i} &= s_i, \\
 \dot{q_i} &= \frac{1}{\tau_0} (\frac{f_i (1-(1-E_0)^{\frac{1}{f_i}})}{E_0} - \frac{q_i}{v^{1-\frac{1}{\alpha}}}), \\
 \dot{v_i} &= \frac{1}{\tau_0} (f_i - v_i^{\frac{1}{\alpha}}), \\
 y_i &= V_0(a_1(1-q_i) - a_2(1-v_i)),
\end{align*}
with the parameters $\{ \tau_s,\tau_f,\tau_0,E_0,\alpha,V_0,a_1,a_2 \}$.
As the parameters $\{ V_0, a_1, a_2 \}$ are not part of the dynamic system they will be excluded from the reduction and estimation and remain fixed at their prior value.

In the scope of this work a linearized fMRI forward sub-model from \cite{kamrani12} is utilized to be applicable in a fully linear setting:
\begin{align*}
 A_{\text{for},i} &= \begin{pmatrix} -\frac{1}{\tau_s} & \frac{1}{\tau_f} & 0 & 0 \\ 1 & 0 & 0 & 0 \\ 0 & \frac{1}{\tau_0 E_0}(1-(1-E_0)(1-\ln(1-E_0))) & -\frac{1}{\tau_0} & \frac{1-\alpha}{\tau_0 \alpha} \\ 0 & \frac{1}{\tau_0} & 0 & -\frac{1}{\tau_0 \alpha} \end{pmatrix}, \\
 B_{\text{for},i} &= \begin{pmatrix} 1 & 0 & 0 & 0 \end{pmatrix}^T, \\
 C_{\text{for},i} &= \begin{pmatrix} 0 & 0 & -a_1 V_0 & a_2 V_0 \end{pmatrix}, \\
 \dot{y_i} &= A_{\text{for},i} + B_{\text{for},i} x_i, \\
 z_i &= C_{\text{for},i} y_i.
\end{align*}
Thus, each state of the common dynamic sub-model, given by a linear control system, has an individual SISO control system attached of which the output reflects fMRI measurements.

The dynamic and the linearized forward sub-models need to be rearranged to fit the linear control system framework:
\begin{align*}
 \dot{X} &= \begin{pmatrix} A_{\text{dyn}} & 0 & 0 & \hdots & 0 \\ \delta_{1,1} & A_{\text{for},1} & 0 \\ \delta_{1,2} & 0 & A_{for,2} \\ \vdots & & & \ddots \\ \delta_{1,n} & & & & A_{\text{for},n} \end{pmatrix} \begin{pmatrix} x \\ y_1 \\ y_2 \\ \vdots \\ y_n \end{pmatrix} + \begin{pmatrix} B_{\text{dyn}} \\ 0 \\ 0 \\ \vdots \\ 0 \end{pmatrix} u, \\
       Z &= \begin{pmatrix} 0 & \begin{pmatrix} C_{\text{for},1} & & \\ & \ddots & \\ & & C_{\text{for},n} \end{pmatrix} \end{pmatrix} \begin{pmatrix} x \\ y_1 \\ \vdots \\ y_n \end{pmatrix},
\end{align*}
with $\delta_{ij}$ being the Kronecker matrix, whose only nonzero element is at ${i,j}$.

For the following experiments the dynamic sub-model's control system is embedded into the fMRI connectivity model.
Since the inference targets the connectivity parameters, each region is assumed to have the same hemodynamic parameters \cite{friston03}.
The number of regions is varied with $n= \{ 9,16,25,36 \}$, which leads to $N = 5n$.
Thus, $P = n^2 + 5$; and as each region is potentially able to receive external input the number of inputs equals the number of regions.
A connectivity matrix $A_{\text{dyn}}$ is generated randomly, but stable and input $u(t)$ will be given by an initial delta impulse.
The prior mean of the parameters is set to $-1$ on the diagonal, $0$ off the diagonal, thus ensuring initial stability of the system ($\theta_{\text{prior}} = -\1$); while the prior covariance is set to the unit matrix $S = \1_{n^2}$.
For the hemodynamic parameters, the prior values assumed for $\{ \tau_s,\tau_f,\tau_0,E_0,\alpha,V_0,a_1,a_2 \}$ are listed in table~\ref{tab:pri}.
In the following applied reduction methods the POD-Greedy state selection will be used, since it seems the most robust for this model.

\begin{table}\centering
\begin{tabular}{c|c|c}
 Parameter & Mean & Covariance \\ \hline
 $\tau_s$ & 0.65 & 0.001 \\
 $\tau_f$ & 0.41 & 0.001 \\
 $\tau_0$ & 0.98 & 0.001 \\
 $E_0$    & 0.34 & 0.001 \\
 $\alpha$ & 0.32 & 0.001 \\
 $V_0$    & 1.00 & 0 \\
 $a_1$    & 1.00 & 0 \\
 $a_2$    & 1.00 & 0 
\end{tabular}
\caption{Prior values for hemodynamic parameters, taken from \cite{friston03}.}
\label{tab:pri}
\end{table}

\begin{figure}[h!]
  \includegraphics[scale=0.38]{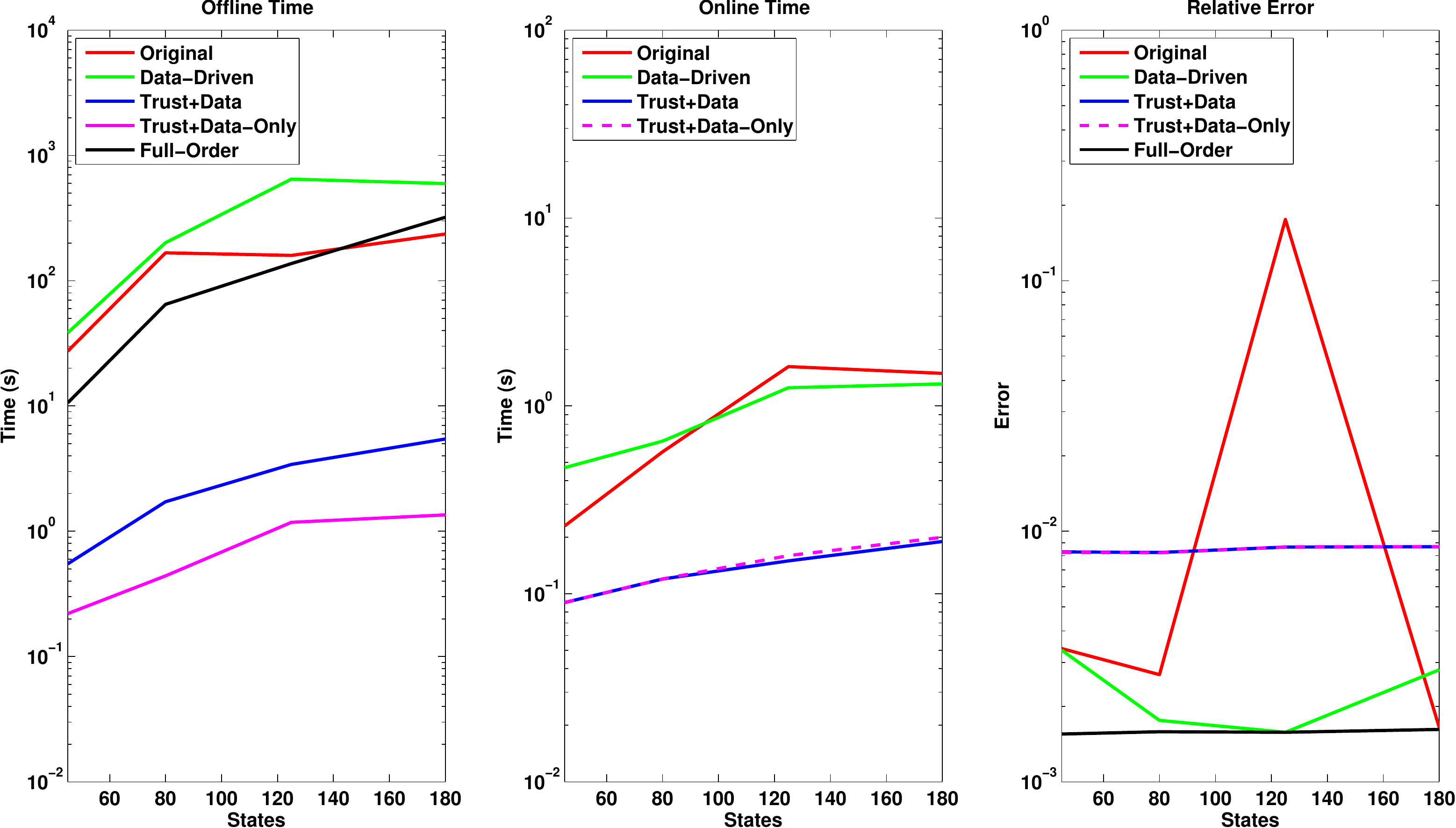}
 \caption{Offline time, online time and relative error for the original, data-driven, trust-region and combined data-driven/trust-region reduction for the fMRI connectivity inference.}
 \label{fig:fmri}
\end{figure}

For all the reductions the number of iterations is fixed to the number of regions, which allows comparison of the different combined reduction variants of the same reduced order.
Figure~\ref{fig:fmri} depicts the results using the same setup as for the previous comparison.

In the offline phase the performance is similar to the generic network, again the original reduction method and data-misfit enhancement visibly consume more time due the relatively higher state dimension.
Here, the data-driven method provides a low relative error, while the original algorithm may produce outliers with high errors.
For these larger models the combination of trust-region and data-driven enhancements results in shorter offline and online phases, yet with higher relative errors.
Again, the exclusive use of data-driven and trust-region ($\alpha=0$) extensions further reduces the offline time without increasing online time or relative error.

Inspecting the total time for solving the inverse problem compared to the full-order and original algorithms' durations is summarized in table~\ref{eff2}. 
Again a speed-up of up to two orders of magnitude is obtained.

\begin{table}\centering
\begin{tabular}{c|c|c}
 State Dimension & Speed-Up (vs. Full-Order) & Speed-Up (vs. Original) \\ \hline
  45 &   $7.42$ &   $3.23$ \\
  80 & $107.50$ &  $10.16$ \\
 125 & $119.43$ &  $33.06$ \\
 180 & $151.81$ & $122.90$
\end{tabular}
\caption{Effectivity of the exclusive combined data-driven and trust-region method for the connectivity model.}
\label{eff2}
\end{table}

\subsection{Extreme-Scale Experiment}
The very short offline phase of the trust-region based reduction process allows a combined reduction of extreme-scale problems \cite{buithanh12}.
As a demonstration of the efficiency under such conditions, we finally look at  the generic model with $N = 256$ states which implies a parameter space dimension of $P = 65536$.
The optimization of the full-order model utilizing an unconstrained optimization, for example a trust-region newton approach\footnote{as MATLAB's \texttt{fminunc}}, 
is not feasible without providing a gradient and Hessian.

In cases with more complicated parametrization, however,  gradients (and Hessians) might not be available or complicated to obtain.
Then, the trust-region-based method presented in this work allows a swift model reduction.
Table~\ref{tab:extr} shows the offline time, online time and relative error in outputs of the reduced order system after inference to the original system.

\begin{table}\centering
 \begin{tabular}{c|c|c|c}
  Method & Offline Time [s] & Online Time [s] & Relative Output Error \\ \hline
  Full Order      & - & - & - \\
  Original        & - & - & - \\
  Data-Driven     & - & - & - \\
  Trust+Data      & 1143.51 & 20.00 & 0.0151 \\
  Trust+Data-Only &  276.74 & 22.31 & 0.0151
 \end{tabular}
 \caption{Offline time, online time and relative error of reduced output using optimized reduced parameters.}
 \label{tab:extr}
\end{table}

The full order optimization as well as original optimization-based reduction method and the data-driven approach were not able to complete the optimization.
This is due to the memory requirements of the unconstrained optimization; since no gradient or Hessian is provided a finite difference scheme is employed to approximate the derivatives.
The required matrices of size $P \times P$ exceeded the test systems memory by far, for the finite difference scheme employed to approximate the derivatives.

As depicted in table \ref{tab:extr},  the combined trust-region/data-driven and trust-region/data-only methods, however, were able to compute a result with $1,5 \%$ relative output error, as never an optimization of the full parameter space had to be performed. Both methods resulted in a comparable online time, while the trust-region/data-driven-only approach has an additional speed-up in the offline-phase, since never a full-order system had to be solved during the Greedy optimizations. We expect a comparable offline performance for the trust-region/data-driven approach, if an efficient to evaluate error estimator would be used as an output error 
measure instead of the true reduction error that was employed in these numerical experiments.

\section{Conclusion and Outlook}
In this contribution we proposed a new data-driven and trust-region approach for projection based model reduction in Bayesian inverse problems. 
Both new ingredients improve the performance of the proposed combined state and parameter model reduction algorithm.
While the data-driven extension is able to reduce relative output errors, the trust-region enhancement massively lowers offline duration in the greedy optimization, 
yet only slightly increases the relative error.
The combination of both results in a shortened offline duration undercutting all previous, but the relative output error corresponds to that of the original reduction algorithm
introduced in \cite{lieberman10}.
Due to the very short offline durations, our new approach allows an efficient inversion of large-scale and even extreme-scale problems as demonstrated in the numerical experiments.

In our numerical experiments, the optimization method inside the reduction algorithms was restricted to an unconstrained optimization due to ensure compatibility for different platforms.
A constrained optimization with provided gradient and Hessian as well as sparse large-scale facilities can additionally improve the optimization during reduction and the actual parameter distribution estimation.

More custom parametrizations than the basic full parametrization can be introduced and thus enable the reduction of more complex models as demonstrated by the fMRI example.
Further research will encompass the generalization of this approach to certain classes of nonlinear models \cite{galbally10}, for example bilinear systems could be reduced with this method requiring only minor adaption.


\section{Acknowledgements}
This work was supported by the Deutsche Forschungsgemeinschaft, DFG EXC 1003 Cells in Motion - Cluster of Excellence, M\"unster, Germany
as well as by the Center for Developing Mathematics in Interaction, DEMAIN, M\"unster, Germany.


\bibliographystyle{unsrt}
\bibliography{optr}

\begin{thebibliography}{10}

\bibitem{friston03}
K.J. Friston, L.~Harrison, and W.~Penny.
\newblock \bibdoi{Dynamic causal modelling}{10.1016/S1053-8119(03)00202-7}.
\newblock {\em Neuroimage}, 19(4):1273--1302, 2003.

\bibitem{stephan07}
K.E. Stephan and K.J. Friston.
\newblock \bibdoi{Models of Effective Connectivity in Neural
  Systems}{10.1007/978-3-540-71512-2_10}.
\newblock In V.K. Jirsa and A.R. McIntosh, editors, {\em Handbook of Brain
  Connectivity}, Understanding Complex Systems, pages 303--327. Springer Berlin
  Heidelberg, 2007.

\bibitem{stuart10}
A.W. Stuart.
\newblock \bibdoi{Inverse problems: a Bayesian
  perspective}{10.1017/S0962492910000061}.
\newblock {\em Acta Numerica}, 19(1):451--559, 2010.

\bibitem{biegler11}
L.~Biegler, G.~Biros, O.~Ghattas, M.~Heinkenschloss, D.~Keyes, B.~Mallick,
  L.~Tenorio, B.~Waanders, K.~Willcox, and Y.~Marzouk.
\newblock {\em \bibdoi{Large-Scale Inverse Problems and Quantification of
  Uncertainty}{10.1002/9780470685853}}.
\newblock Wiley Series in Computational Statistics. Wiley, 2011.

\bibitem{WMGGD99}
S.~Weile, E.~Michielssen, E.~Grimme, and K.~Gallivan.
\newblock \bibdoi{A method for generating rational interpolant reduced order
  models of two- parameter linear systems}{10.1016/S0893-9659(99)00063-4}.
\newblock {\em Appl. Math. Lett.}, 12(5):93--102, 1999.

\bibitem{FB07}
L.~Feng and P.~Benner.
\newblock \bibdoi{A Robust Algorithm for Parametric Model Order
  Reduction}{10.1002/pamm.200700749}.
\newblock In {\em PAMM}, volume 7(1), page 1021501–1021502, 2007.

\bibitem{BB08}
U.~Baur and P.~Benner.
\newblock \biburl{Parametrische Modellreduktion mit d{\"u}nnen
  Gittern}{http://www.de.mpi-magdeburg.mpg.de/mpcsc/benner/pub/BaurBenner-GMA-Proceedings2008.pdf}.
\newblock In {\em GMA-Fachausschuss 1.30, Modellbildung, Identifizierung und
  Simulation in der Automatisierungstechnik, Salzburg ISBN 978-3-9502451-3-4},
  2008.

\bibitem{LE09}
B.~Lohmann and R.~Eid.
\newblock \biburl{Efficient Order Reduction of Parametric and Nonlinear Models
  by Superposition of Locally Reduced
  Models}{http://www.rt.mw.tum.de/fileadmin/w00bhf/www/publikationen/2009_Lohmann_Hirschberg.pdf}.
\newblock In {\em Methoden und Anwendungen der Regelungstechnik.
  Erlangen-M{\"u}nchener Workshops 2007 und 2008}. Shaker Verlag, Aachen, 2009.

\bibitem{eid11}
R.~Eid, R.~Casta{\~n}{\'e}-Selga, H.~Panzer, T.~Wolf, and B.~Lohmann.
\newblock \bibdoi{Stability-preserving parametric model reduction by matrix
  interpolation}{10.1080/13873954.2011.547671}.
\newblock {\em MATH. COMP. MODEL. DYN.}, 17(4):319--335, 2011.

\bibitem{amsallem11}
D.~Amsallem and C.~Farhat.
\newblock \bibdoi{An online method for interpolating linear parametric
  reduced-order models}{10.1137/100813051}.
\newblock {\em SIAM J. Sci. Comput.}, 33(5):2169--2198, 2011.

\bibitem{haasdonk08}
B.~Haasdonk and M.~Ohlberger.
\newblock \bibdoi{Reduced basis method for finite volume approximations of
  parametrized linear evolution equations}{10.1051/m2an:2008001}.
\newblock {\em M2AN}, 42(2):277--302, 2008.

\bibitem{haasdonk13}
B.~Haasdonk.
\newblock \bibdoi{Convergence rates of the {POD}-greedy
  method}{10.1051/m2an/2012045}.
\newblock {\em M2AN}, 47(3):859--873, 2013.

\bibitem{haasdonk11}
B.~Haasdonk and M.~Ohlberger.
\newblock \bibdoi{Efficient reduced models and {\it a posteriori} error
  estimation for parametrized dynamical systems by offline/online
  decomposition}{10.1080/13873954.2010.514703}.
\newblock {\em MATH. COMP. MODEL. DYN.}, 17(2):145--161, 2011.

\bibitem{nguyen11}
N.C. Nguyen, G.~Rozza, D.B.P. Huynh, and A.T. Patera.
\newblock \biburl{Reduced basis approximation and a posteriori error estimation
  for parametrized parabolic {PDE}s: application to real-time {B}ayesian
  parameter
  estimation}{http://augustine.mit.edu/methodology/papers/atpWileyMay2009preprint.pdf}.
\newblock In {\em Large-scale inverse problems and quantification of
  uncertainty}, Wiley Ser. Comput. Stat., pages 151--177. Wiley, 2011.

\bibitem{himpe13}
C.~Himpe and M.~Ohlberger.
\newblock \biburl{Cross-Gramian Based Combined State and Parameter Reduction
  for Large-Scale Control Systems}{http://arxiv.org/abs/1302.0634}.
\newblock arxiv (math.oc) 1302.0634, preprint (submitted), Institute for
  Computational and Applied Mathematics, 2013.

\bibitem{lieberman10}
C.~Lieberman, K.~Willcox, and O.~Ghattas.
\newblock \bibdoi{Parameter and state model reduction for large-scale
  statistical inverse problems}{10.1137/090775622}.
\newblock {\em SIAM J. Sci. Comput.}, 32(5):2523--2542, 2010.

\bibitem{antoulas09}
A.C. Antoulas.
\newblock \bibdoi{An overview of model reduction methods and a new
  result}{10.1109/CDC.2009.5400920}.
\newblock In {\em Proceedings of the 48th IEEE Conference on Decision and
  Control}, pages 5357--5361. IEEE, 2009.

\bibitem{buithanh07}
T.~Bui-Thanh, K.~Willcox, O.~Ghattas, and B.~van Bloemen~Waanders.
\newblock \bibdoi{Goal-oriented, model-constrained optimization for reduction
  of large-scale systems}{10.1016/j.jcp.2006.10.026}.
\newblock {\em J. Comput. Phys.}, 224(2):880--896, 2007.

\bibitem{buithanh08}
T.~Bui-Thanh, K.~Willcox, and O.~Ghattas.
\newblock \bibdoi{Model reduction for large-scale systems with high-dimensional
  parametric input space}{10.1137/070694855}.
\newblock {\em SIAM J. Sci. Comput.}, 30(6):3270--3288, 2008.

\bibitem{bashir08}
O.~Bashir, K.~Willcox, O.~Ghattas, B.~van Bloemen~Waanders, and J.~Hill.
\newblock \bibdoi{Hessian-based model reduction for large-scale systems with
  initial-condition inputs}{10.1002/nme.2100}.
\newblock {\em Int. J. Numer. Meth. Engng.}, 73(6):844--868, 2008.

\bibitem{lieberman12}
C.~Lieberman and K.~Willcox.
\newblock \bibdoi{Goal-Oriented Inference: Approach, Linear Theory, and
  Application to Advection Diffusion}{10.1137/110857763}.
\newblock {\em SIAM J. Sci. Comput.}, 34(4):A1880--A1904, 2012.

\bibitem{arian00}
E.~Arian, M.~Fahl, and E.W. Sachs.
\newblock \biburl{Trust-region proper orthogonal decomposition for flow
  control}{http://www.dtic.mil/dtic/tr/fulltext/u2/a377382.pdf}.
\newblock Technical report, DTIC Document, 2000.

\bibitem{lieberman07}
C.~Lieberman and B.~Van Bloemen~Waanders.
\newblock \biburl{Hessian-Based Model Reduction Approach to Solving Large-Scale
  Source Inversion Problems}{http://csri.sandia.gov/Proceedings/CSRI2007.pdf}.
\newblock In {\em CSRI Summer Proceedings 2007}, pages 37--48, 2007.

\bibitem{grepl11}
Martin~A. Grepl and Mark Kärcher.
\newblock \bibdoi{Reduced basis a posteriori error bounds for parametrized
  linear-quadratic elliptic optimal control
  problems}{http://dx.doi.org/10.1016/j.crma.2011.07.010}.
\newblock {\em Comptes Rendus Mathematique}, 349(15-16):873--877, 2011.

\bibitem{lall99}
S.~Lall, J.E. Marsden, and S.~Glavaski.
\newblock \biburl{Empirical model reduction of controlled nonlinear
  systems}{authors.library.caltech.edu/20343/2/10.1.1.123.4669.pdf}.
\newblock {\em Proceedings of the IFAC World Congress}, F:473--478, 1999.

\bibitem{lall02}
S.~Lall, J.E. Marsden, and S.~Glavaski.
\newblock \bibdoi{A subspace approach to balanced truncation for model
  reduction of nonlinear control systems}{10.1002/rnc.657}.
\newblock {\em INT. J. ROBUST. NONLIN.}, 12(6):519--535, 2002.

\bibitem{kelley99}
C.T. Kelley and E.W. Sachs.
\newblock \bibdoi{A trust region method for parabolic boundary control
  problems}{10.1137/S1052623496308965}.
\newblock {\em SIAM J. OPTIMIZ.}, 9(4):1064--1081, 1999.

\bibitem{conn00}
A.R. Conn, N.I.M. Gould, and P.L. Toint.
\newblock {\em \biburl{Trust Region
  Methods}{http://books.google.de/books?id=5kNC4fqssYQC}}.
\newblock MPS-SIAM Series on Optimization. SIAM, 2000.

\bibitem{kamrani12}
E.~Kamrani, A.~Foroushani, M.~Vaziripour, and M.~Sawan.
\newblock \bibdoi{Detecting the Stable, Observable and Controllable States of
  the Human Brain Dynamics}{10.4236/ojmi.2012.24024}.
\newblock {\em OJMI}, 2(4):128--136, 2012.

\bibitem{buithanh12}
T.~Bui-Thanh, C.~Burstedde, O.~Ghattas, J.~Martin, G.~Stadler, and L.C. Wilcox.
\newblock \bibdoi{Extreme-scale UQ for Bayesian inverse problems governed by
  PDEs}{10.1137/110845598}.
\newblock {\em Proceedings of the International Conference on High Performance
  Computing, Networking, Storage and Analysis}, page~3, 2012.

\bibitem{galbally10}
D.~Galbally, K.~Fidkowski, K.~Willcox, and O.~Ghattas.
\newblock \bibdoi{Non-linear model reduction for uncertainty quantification in
  large-scale inverse problems}{10.1002/nme.2746}.
\newblock {\em Int. J. Numer. Meth. Engng.}, 81(12):1581--1608, 2010.

\end{thebibliography}

\end{document}